\documentclass{amsart}
\usepackage{latexsym,amssymb,graphics}

\def\setof#1#2{\{\,{#1}\,:\, {#2}\,\}}

\newtheorem{theorem}{Theorem}

\newtheorem{lemma}{Lemma}

\newtheorem{corollary}{Corollary}
\newtheorem{example}{Example}

\begin{document}

\Large

\title[Monotonic MLE]{Maximum-Likelihood Non-Decreasing Response Estimates}
\author {L. Thomas Ramsey}
\address{Department of Mathematics\\
University of Hawaii at Manoa \\
Honolulu, HI 96822}
\email{ramsey@math.hawaii.edu}
\urladdr{http://www.math.hawaii.edu/~ramsey}
\thanks{Mahalo nui loa to David Z. Liang.}
\date{July 5, 2011}

\begin{abstract}  Let $x_{i,j}$, $1 \le i \le m$, $1 \le j \le n_i$, be observations from a doubly-indexed sequence $\{X_{i,j}\}$ of independent random variables (all of them discrete, or all of them absolutely continuous).  Suppose that each $X_{i,j}$ has the PDF $f(x\mid\theta_i)$
from a one-parameter family of PDFs $f(x\mid \theta)$.  Mild assumptions are described under which
there is a unique compound estimate $\mathbf \phi=\langle \phi_1, \hdots \phi_m \rangle$ of $\mathbf \theta=\langle \theta_1, \hdots
\theta_m \rangle$  such that
\begin{enumerate}
\item
For integers $i<k$ in $[1,m]$, $\phi_i \le \phi_k$ ($\phi$ is non-decreasing with respect to the index $i$).
\item
Among all non-decreasing vectors of parameters $\lambda=\langle \lambda_1, \hdots, \lambda_m\rangle$,
$$\lambda \ne
\phi \quad \Rightarrow \quad \ell(\mathbf x\mid\lambda) < \ell(\mathbf x\mid\phi)
$$
 where
\begin{enumerate}
\item
$\mathbf x$ is the (doubly-indexed) vector of observations
\item
$\ell$ is the compound likelihood function: 
$$
\ell(\lambda\mid\mathbf x)=\prod_{i=1}^m \prod_{j=1}^{n_i} f(x_{i,j}\mid\lambda_i)
$$
\end{enumerate}
\end{enumerate}
An efficient algorithm is described to compute $\mathbf \phi$.  The notation $[a \hdots b]$ denotes
the integers in the real interval $[a,b]$.  For $J  \subset [1 \hdots m]$ let
$$
\mu\,J=\mu(J)=\dfrac{\sum_{i\in J}\sum_{j=1}^{n_i} x_{i,j}}{\sum_{i=a}^b n_i}
$$
That is, $\mu(J)$ is the sample mean of observations $x_{i,j}$ with $i \in J$.  Here is the theorem that justifies the algorithm:
\begin{enumerate}
\item
Let $\tau_1<\tau_2< \hdots \tau_s$ be a complete list of the distinct components
of $\phi$.  Let $A_r=\setof{i \in [1 \hdots m]}{\phi_i=\tau_r}$. 
There are integers $a_r \le b_r$ such that
$A_r=[a_r \hdots b_r]$.  Also, $\tau_r=\mu(A_r)$. 
\item
For integers $r \in [1 \hdots s]$,
set
$$
\kappa_r=\min \setof{\mu[a_r\, \hdots \,k]}{k \in [a_r\,  \hdots \,m]}
$$
Then,
$$
b_r=\max\setof{k \in [\,a_r\, \hdots \,m]}{\mu[\,a_r\, \hdots\, k]=\kappa_r}
$$
\end{enumerate}
\end{abstract}

\maketitle

\section{Introduction and Formal Context}

 Let $x_{i,j}$, $1 \le i \le m$, $1 \le j \le n_i$, be observations from a doubly-indexed sequence $\{X_{i,j}\}$ of independent random real variables.  Suppose that each $X_{i,j}$ has the PDF $f(x\mid\theta_i)$
from a one-parameter family of PDFs $\setof{f(x\mid \theta)}{\theta \in \Theta}$.

In this note, we make {\bf compound} maximum likelihood estimations $\widehat \theta$ of $\theta=\langle \theta_1, \hdots \theta_m\rangle$, 
subject only to a non-decreasing constraint:
$$
i \le k \quad \Rightarrow \quad \widehat {\theta_i} \le \widehat{\theta_k}
$$
Compound estimates $\widehat \theta \in \Theta^m$ that meet this constraint will be called {\bf non-decreasing}.

Think of each value of the index $i$ specifiying consecutive levels of an explanatory variable, and $\theta_i$ is the response to that $i$-th level of the explanatory
variable.  The goal of this note is to specify the maximum likelihood non-decreasing response estimate (existence and uniqueness) and provide an algorithm for its efficient computation.

Of course, the same technology gives also the maximum likelihood non-increasing response function.  One simply reverses the ordering of the levels of
the explanatory variables, and applies the same theory and algorithm to the reversely ordered data.

\subsection{Formal Context (Assumptions)}  Let $x_{i,j}$, $1 \le i \le m$, $1 \le j \le n_i$, be observations from a doubly-indexed sequence $\{X_{i,j}\}$ of independent random real variables.  Suppose that each $X_{i,j}$ has the PDF $f(x\mid\theta_i)$
from a one-parameter family $\mathcal H$ of PDFs $\setof{f(x\mid \theta)}{\theta \in \Theta}$.  Set
$$
\mathcal D=\setof{x \in \mathbb R}{(\exists \theta \in \Theta)(f(x \mid \theta)>0)}
$$
A real number $x$ is called {\bf observable} if and only if $x \in \mathcal D$.

{\bf Assumption 1.}  One of two cases are assumed:
\begin{enumerate}
\item
The PDFs in $\mathcal H$ are for discrete real random variables.
\item
The PDFs in $\mathcal H$ are for absolutely continuous real random variables.
\end{enumerate}

{\bf Assumption 2.}  We assume that $\mathcal D \subset \Theta$ and $\Theta$ is a real interval of positive length.  Consequently, $\Theta$ includes the arithmetic
means of finite sequences of observable real numbers.

{\bf Assumption 3.}  Let $y_s$, $1 \le s \le t$ be observations from the independent random variables $Y_s$, $1 \le s \le t$, with PDFs in the given family.  As in \cite{statisticalinference2002} and on pages 337 -- 341 of \cite{lossmodels2004}, the likelihood function is
$$
L(\theta|\mathbf y)= \prod_{s=1}^t f(y_i\mid \theta)
$$
where $\mathbf y =\langle y_1 \hdots y_t \rangle$.  We assume that, if each $y_i$ is observable and
$
\overline y
$ is the arithmetic mean of $y_1$, $\hdots$, $y_t$, then
 $L(\theta \mid \mathbf y)$ is strictly increasing for  $\theta\le
\overline y$ and strictly decreasing for $\theta \ge \overline y$.  

The next lemma is an immediate consequence of these assumptions.  It will be applied
frequently to subsets of the observations under discussion.
\begin{lemma}\label{L:begin}   Let $y_s \in \mathcal D$, $1 \le s \le t$, be observations from the independent random variables $Y_s$, $1 \le s \le t$, with PDFs in the given family. Then  $\overline y \in \Theta$ and,
for all $\eta \in \Theta$ such that $\eta \ne \overline y$
\begin{equation}
L(\overline y \mid \mathbf y)>L(\eta \mid \mathbf y)
\label{E:eqn0}
\end{equation}
Also,
$ L(\overline y \mid \mathbf y)>0.$
\end{lemma}
\begin{proof}  by Assumption 2, $\overline y \in \Theta$.  Equation \ref{E:eqn0} follows directly from Assumption 3.  Since
$\Theta$ is a real interval of positive length, there is some $\eta \ne \overline y$ in
$\Theta$.  Therefore
$$
L(\overline y \mid \mathbf y) >L(\eta \mid \mathbf y) \ge 0
$$
\end{proof}

\subsection{Applicability of the Formal Context}  The formal context applies to many of the common one-parameter distributions, if one parameterizes them by their means and takes care to include some boundary distributions.  Here are some examples, verified in an appendix:
\begin{enumerate}
\item
Bernoulli random variables with parameter $p \in [0,1]$.  Note that we include the boundary cases $0$ and $1$ because they could be sample means.   This follows the usual practice for
maximum likelihood estimators (see page 318 of \cite{statisticalinference2002}).
\item
Poisson random variables parameterized by their means, with $\Theta=[0,\infty)$ (see page 644 of \cite{lossmodels2004}).  Note the inclusion of the boundary value $0$.
\item
Let $\Theta=[0, \infty)$ and $\mathcal F$ the family of geometric random variables (here including one constant random variable), parameterized by
their means (see page 644 of \cite{lossmodels2004}).
Given $\theta \in \Theta$, set $p=1/(1+\theta)$.  For $\theta >0$, the PDF with parameter $\theta$ is defined as follows:
$$
f(x\mid\theta) =\begin{cases}
(1-p)^x p, & \text{if $x \ge 0$ and an integer} \\
\\
0 &\text{otherwise}
\end{cases}
$$
For $\theta=0$, $f(x \mid 0)=0$ for all $x$ except that $f(0 \mid 0)=1$.
\item
Normal distributions $N(\theta,\sigma)$ with $\sigma>0$ fixed and $\theta \in \mathbb R$.  It is noted on page 317 of \cite{statisticalinference2002} that $L(\theta \mid \mathbf y)$
has a unique maxium at $\theta =\overline y$.  To verify the slightly stronger assumptions of this note, one slightly modifies the argument given in \cite{statisticalinference2002}.
\item
Exponential random variables parameterized by their means (page 638 of \cite{lossmodels2004}).
Let $\Theta=(0, \infty)$.
Given $\theta \in \Theta$, the PDF with parameter $\theta$ is defined as follows:
$$
f(x\mid\theta) =\begin{cases}
\dfrac 1 \theta  \exp(-x/\theta), & \text{for  $x > 0$} \\
\\
0 &\text{for $x \le 0$}
\end{cases}
$$
\end{enumerate}

The motivating example of the next section has binomial distributions.  The example is an application of a corresponding Bernoulli random variable case.  The compound likelihood functions of the two situations differ by a non-zero factor $C$ that does not depend on any of the parameters $p_i$.  $C$ has the form
$$
\prod_{i=1}^m \binom{n_i}{d_i}
$$
where $d_i$ is the number of successes ($1$s) in the Bernoulli observations $x_{i,1}$, $\hdots$, $x_{i,n_i}$.

\subsection{The Algorithm}  Within this formal context, and with each $x_{i,j}$ observable, it will be proved that there is a unique non-decreasing compound estimate $\mathbf \phi=\langle \phi_1, \hdots \phi_m \rangle$ of $\mathbf \theta=\langle \theta_1, \hdots
\theta_m \rangle$  such that, among all non-decreasing  $\lambda\in \Theta^m$,
$$\lambda \ne
\phi \quad \Rightarrow \quad \ell(\mathbf x\mid\lambda) < \ell(\mathbf x\mid\phi)
$$
 where
\begin{enumerate}
\item
$\mathbf x$ is the (doubly-indexed) vector of observations
\item
$\ell$ is the compound likelihood function: 
$$
\ell(\lambda\mid\mathbf x)=\prod_{i=1}^m \prod_{j=1}^{n_i} f(x_{i,j}\mid\lambda_i)
$$
\end{enumerate}
There is an efficient algorithm to compute $\mathbf \phi$.  The notation $[a \hdots b]$ denotes
the integers in the $[a,b]$.  For $J  \subset [1 \hdots m]$ let
$$
\mu\,J=\mu(J)=\dfrac{\sum_{i\in J}\sum_{j=1}^{n_i} x_{i,j}}{\sum_{i=a}^b n_i}
$$
That is, $\mu(J)$ is the sample mean of observations $x_{i,j}$ with $i \in J$.

Let $\tau_1<\tau_2< \hdots \tau_s$ be a complete list of the distinct components
of $\phi$.  Let $A_r=\setof{i \in [1 \hdots m]}{\phi_i=\tau_r}$. 
By a theorem, there are integers $a_r \le b_r$ such that
$A_r=[a_r \hdots b_r]$.  Also, $\tau_r=\mu(A_r)$.

For integers $r \in [1 \hdots s]$,
set
\begin{equation}
\kappa_r=\min \setof{\mu[a_r\, \hdots \,k]}{k \in [a_r\,  \hdots \,m]}
\label{E:kappa}
\end{equation}
It will be proved that
\begin{equation}
b_r=\max\setof{k \in [\,a_r\, \hdots \,m]}{\mu[\,a_r\, \hdots\, k]=\kappa_r}
\label{E:br}
\end{equation}

Here is the algorithm:
\begin{enumerate}
\item
Set $a_1$ equal to $1$.  Compute $\kappa_1$ according to Equation \ref{E:kappa}.  Then compute $b_1$ according to Equation \ref{E:br}.
For $i \in [a_1 \hdots b_1]$, set $\phi_i$ equal to $\mu[a_1 \hdots b_1]$.
\item
Proceed recursively.  Suppose that $\{a_t\}_{t=1}^r$ and $\{b_t \}_{t=1}^r$ satisfy Equation $\ref{E:br}$ with $a_1=1$, and 
$a_{t}=b_{t-1}+1$ for $1<t \le r$.  If $b_r=m$, the algorithm stops.  If $b_r<m$, set $a_{r+1}=b_r+1$.  Then compute $\kappa_{r+1}$ according to Equation $\ref{E:kappa}$.  Then compute $b_{r+1}$ according to Equation \ref{E:br}.  For $i \in [a_{r+1},b_{r+1}]$, set $\phi_i$ equal to $\mu[a_{r+1} \hdots
b_{r+1}]$.
\end{enumerate}
\newpage

\section{A Motivating Example}

For a particular mathematics course, can one use various SAT scores to predict performance on the final examination?  A significant subset of the students failed
to take the final examination, in some cases because they withdrew from the course (with a record of W) or they simply received a letter grade of F.
In the Data Appendix, Table \ref{tab:data} tabulates the no-show counts by SAT-R levels.  The data are noisy. 

\begin{figure}[htb!]
\centering%
\includegraphics{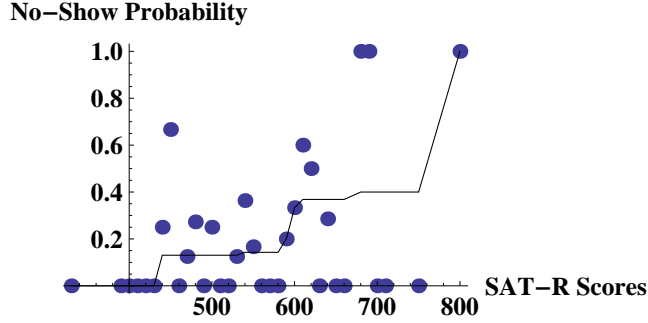}
\caption{ML Non-Decreasing Response Estimate}
\label{fig:NoShowVsRSAT}
\end{figure}

Table \ref{tab:noshownd} gives the maximum likelihood non-decreasing
response estimate.  In Figure  \ref{fig:NoShowVsRSAT} are plotted both the observed no-show rates for each SAT-R score, and the maximum-likelihood non-decreasing response estimate.

On the other hand, the maximum likelihood {\bf non-increasing} response estimate is the constant $26/152$ for every SAT-R level.

To estimate the significance of the sharp difference between the maximum-likelihood non-decreasing and non-increasing response estimates,
10,000 data tables with the same structure as Table \ref{tab:data} were generated under a null hypothesis of a constant no-show rate of
$26/152$ regardless of SAT-R score.
\begin{itemize}
\item
For each simulated table, the maximum-likelihood non-decreasing response estimate $f$ and the maximum-likelihood non-increasing response estimate $g$ were computed.
This test statistic was tabulated:
$$
\Delta=[ \, f(800)-f(330) \,] - [g(330)-g(800)]
$$
\item
The simulated quantile ranking is   0.9899  for the same statistic for the actual table (namely $1$).  That is, 101 of the simulated tables had $\Delta=1$ (the maximum
logically possible).
\end{itemize}

There are easier and more routine ways to reject the null hypothesis, but the non-decreasing response estimate itself has proved useful for directing resources to students who have a greater risk of failure.  

The observed table is broadly consistent with an alternative hypothesis of the maximum-likelihood non-decreasing response estimate; its simulated loglikelihood rank (based on 10,000 simulated tables generated under this alternative hypothesis)  is 0.3476.

\begin{table}
\centering
\begin{tabular}{cccc}
\text{SAT-R Scores} &\text{Total Count} &\text{No-Show Count} &\text{No-Show Rate} \\
330-430 &15 & 0 & 0\% \\
440-530 &69 &9 & 13.0\% \\
540-580 &35 & 5 &14.3\% \\
590 &5 & 1 & 20\% \\
600 & 3 & 1 & 33.3\%\\
610-660 & 19 & 7 &36.8\% \\
680-750 &5 & 2 &40\% \\
800 & 1 & 1 & 100\%\\
\\
\text{Overall:} &152 &26 & 17.1\% 
\end{tabular}
\caption{ML Non-Decreasing Response Estimate of No-Show Rates}
\label{tab:noshownd}
\end{table}

\newpage

\section{Existence of Likelihood-Maximizing Non-Decreasing Response Estimates}

Every $\theta \in \Theta^m$ defines a partition $S(\theta)$ of $[1 \hdots m]$ as follows.  Let $\tau_1$, $\tau_2$, $ \hdots$, $\tau_s$ be a  listing of the distinct
components of $\theta$ without the repetition of any component.  Let
$$
S(\theta)=\{A_1, \hdots, A_s \} \quad \text{where} \quad A_t=\setof{i \in [1 \hdots m]}{ \theta_i = \tau_t }
$$
Let  $I(\theta)$, the index of $\theta$, be equal to the number of $t \in [1 \hdots s]$ such that
$$
\tau_t \ne \mu(A_t)
$$

\begin{lemma} \label{L:step0} Let $\mathbf x=\{\{x_{i,j}\}_{j=1}^{n_i}\}_{i=1}^m$ have each $x_{i,j} \in \mathcal D$. Suppose
that  $\theta \in \Theta^m$ is non-decreasing and has $I(\theta)>0$.
Then there is some non-decreasing $\widetilde \theta  \in \Theta^m$ such that
\begin{enumerate}
\item
$S(\widetilde \theta)$ has fewer members than $S(\theta)$ or, $S(\theta)=S(\widetilde \theta)$ and $I(\widetilde \theta)<I(\theta)$.
\item
For the compound likelihood function $\ell(\cdot \mid \mathbf x)$, 
$$
\ell(\theta \mid \mathbf x) < \ell(\widetilde{\theta} \mid \mathbf x)
$$
\end{enumerate}
\end{lemma}
\begin{proof}  Suppose first that $\ell(\theta \mid \mathbf x)=0$.  Let $\overline x=\mu[1 \hdots m]$, the arithmetic mean of all the observations.  By Lemma \ref{L:begin}, because every $x_{i,j}$ is observable, we have $\overline x \in \Theta$.
Let $\widetilde \theta \in \Theta^m$ be the constant vector with components equal
to $\overline x$.  Clearly $\widetilde \theta$ is non-decreasing, $S(\widetilde \theta)=\{[1 \hdots m]\}$, and for all $i$ we have 
$$
\widetilde{\theta_i}=\overline x=\mu[1 \hdots m]
$$
Thus $I(\theta)=0$.  By Lemma \ref{L:begin},
$$
\ell(\widetilde \theta \mid \mathbf x)=L(\overline x \mid \mathbf x)>0=\ell(\theta \mid \mathbf x)
$$

For the rest of the proof we assume that $\ell(\theta \mid \mathbf x)>0$.

Let $\tau_1<\tau_2 \hdots  <\tau_s$ be a complete list of the distinct components of $\theta$.  Set
$$
A_t=\setof{i \in [1 \hdots m]}{\theta_i=\tau_t}
$$
Then $S(\theta)=\{A_1, \hdots, A_s\}$.  

Because $I(\theta)>0$, there is at least one $T \in [1 \hdots s]$ such that $\tau_T \ne \mu(A_T)$.

Because $\theta$ is non-decreasing, there are integers $a \le b$ in $[1 \hdots m]$ such that $A_T=[a \hdots b]$.

 Let $\mathbf y$ be the vector with $n_a + \hdots +n_b$ components $\{\{x_{i,j}\}_{j=1}^{n_i}\}_{i=a}^b$.  Then,
$$
\ell(\theta \mid \mathbf x)=L(\tau_T \mid \mathbf y) \cdot C
$$
where
$$
L(w \mid \mathbf y)= \prod_{i=a}^b \prod_{j=1}^{n_i} f(x_{i,j} \mid w)
$$
and
$$
C=\prod_{t \in [1 \hdots s]\backslash\{T\}} \prod_{i \in A_t} \prod_{j=1}^{n_i} f(x_{i,j} \mid \theta_i)
$$

Since $\ell(\theta \mid \mathbf x)>0$ and $C \ge 0$, we have $C>0$.

For $t \in [1 \hdots s]\backslash \{T\}$ and for all $i \in A_t$,
we set
$$\widetilde {\theta_i}=\tau_t=\theta_i$$

For $t=T$, we will select some $\widetilde{\tau_T}$ and for all $i \in A_T$ we will set $\widetilde{\theta_i}=\widetilde{\tau_T}$.
Note that we will then have
$$
\ell(\widetilde \theta \mid \mathbf x)=L(\widetilde{\tau_T} \mid \mathbf y) \cdot C
$$
We break the selection of $\widetilde{\tau_T}$ into 6 cases:
\begin{enumerate}
\item
$\widetilde{\tau_T}=\mu(A_1)$ if $s=1$.
\item
$\widetilde{\tau_T}=\mu(A_1)$ if $s>1$, $T=1$ and $\mu(A_1)<\tau_2$.
\item
$\widetilde{\tau_T}=\mu(A_T)$ if $s>1$, $T=s$, and $\mu(A_s)>\tau_{s-1}$.
\item
$\widetilde{\tau_T}=\tau_{T+1}$ if $s>1$, $1 \le T <s$ and $\mu(A_{T})\ge \tau_{T+1}$.
\item
$\widetilde{\tau_T}=\tau_{T-1}$ if $s>1$, $1<T\le s$, and $\mu(A_T) \le \tau_{T-1}$.
\item
$\widetilde{\tau_T}=\mu(A_T)$ if $s>1$, $1<T<s$, and
$$
\tau_{T-1} <\mu(A_T) <\tau_{T+1}
$$
\end{enumerate}
Note that in every case, $\widetilde{\tau_T} \in \Theta$ and thus $\widetilde \theta \in
\Theta^m$:
\begin{itemize}
\item
In cases 1, 2, 3 and 6, $\widetilde{\tau_T}$ is the arithmetic mean of observable real numbers.  By Lemma \ref{L:begin}, $\widetilde{\tau_T}$ is in $\Theta$.
\item
In cases 4 and 5, $\widetilde{\tau_T}$ is equal to a component of $\theta$.
Since $\theta \in \Theta^m$, we have  $\widetilde{\tau_T}$ in $\Theta$.
\end{itemize}

Next, please note that in every case $\widetilde{\tau_T}$ has been carefully selected
to make $\widetilde \theta$ be non-decreasing.

It will be argued shortly that the change from $\theta_i=\tau_T$ to $\widetilde {\theta_i}=\widetilde{\tau_T}$ for $i \in A_T$ 
strictly increases $L(\cdot \mid \mathbf y)$.  Since $C>0$, that will give us
$$
\ell(\widetilde \theta \mid \mathbf x)=L(\widetilde{\tau_T} \mid \mathbf y) \cdot C
>L({\tau_T} \mid \mathbf y) \cdot C=\ell(\theta \mid \mathbf x)
$$
\begin{itemize}
\item
In cases 1, 2, 3 and 6, the change is from $\tau_T \ne \mu(A_T)$ to $\mu(A_T)$.
By Lemma \ref{L:begin},
$$
L(\mu(A_T) \mid \mathbf y)>L(\tau_T \mid \mathbf y)
$$
\item
In case 4, we have $\tau_T <\tau_{T+1} \le \mu(A_T)$.  By Assumption 3,
$L(w \mid \mathbf y)$ is strictly increasing for $w \le \mu(A_T)$.  Therefore
$$
L(\tau_{T+1} \mid \mathbf y) > L(\tau_T \mid \mathbf y)
$$
In case 4, $\widetilde{\tau_T}=\tau_{T+1}$ and thus
$$
L(\widetilde{\tau_{T} }\mid \mathbf y) > L(\tau_T \mid \mathbf y)
$$
\item
In case 5, we have $\mu(A_T) \le \tau_{T-1} <\tau_T$.  By Assumption 3,
$L(w \mid \mathbf y)$ is strictly decreasing for $w \ge \mu(A_T)$.  Therefore,
$$
L(\tau_{T-1} \mid \mathbf y) > L(\tau_T \mid \mathbf y)
$$
In case 5, $\widetilde{\tau_T}=\tau_{T-1}$ and thus
$$
L(\widetilde{\tau_{T}} \mid \mathbf y) > L(\tau_T \mid \mathbf y)
$$
\end{itemize}

Lastly, consider Item (i) of the Lemma.
\begin{itemize}
\item
In case 1, $I(\widetilde \theta)=0$.  Also $S(\widetilde \theta)$ has one member.
If $S(\theta)$ has more than one member, Item (i) is satisfied.  If $S(\theta)$ has one
member, then $S(\theta)=S(\widetilde \theta)=\{[1 \hdots m]\}$.  Since $I(\theta)>0$ by hypothesis, we have
$I(\widetilde \theta)<I(\theta)$.  That, too, satisfies Item (i).
\item
In cases 2, 3 and 6, $S(\theta)=S(\widetilde \theta)$.  Also, for $r \ne T$, there is no
change in whether $\tau_r$ is equal to $\mu(A_r)$.  However, for $i \in A_T$,
by changing
for $\theta_i=\tau_T \ne \mu(A_T)$ to $\widetilde{\theta_i}=\mu(A_T)$,
we have made $I(\widetilde \theta)=I(\theta)-1$.
\item
In case 4, $S(\widetilde \theta)$ has one fewer partition member than $S(\theta)$,
because
$$
S(\widetilde \theta)=\left(S(\theta)\backslash\{A_{T},A_{T+1}\} \right) \cup \{\widetilde{A_{T}}\}
$$
where
$$\widetilde{A_{T}}=A_{T} \cup A_{T+1}
$$
\item
In case 5, again $S(\widetilde \theta)$ has one fewer partition member than $S(\theta)$,
because
$$
S(\widetilde \theta)=\left(S(\theta)\backslash\{A_{T-1},A_T\} \right) \cup \{\widetilde{A_{T-1}}\}
$$
where
$$\widetilde{A_{T-1}}=A_{T-1} \cup A_T
$$
\end{itemize}

\end{proof}

\begin{lemma} \label{L:step1} Let $\mathbf x=\{\{x_{i,j}\}_{j=1}^{n_i}\}_{i=1}^m$ have each $x_{i,j} \in \mathcal D$. Suppose
that  $\theta \in \Theta^m$ is non-decreasing with $I(\theta)>0$.
Then there is some non-decreasing $\widetilde \theta  \in \Theta^m$ with $I(\theta)=0$ such that
$$
\ell(\theta \mid \mathbf x) < \ell(\widetilde{\theta} \mid \mathbf x)
$$
\end{lemma}
\begin{proof}  Apply Lemma \ref{L:step0} to $\theta^{(0)}=\theta$ to get $\theta^{(1)}$, and continue recursively to apply Lemma \ref{L:step0} to 
$\theta^{(k)}$ to get $\theta^{(k+1)}$, but stop if $I(\theta^{(k)})=0$.  

In Lemma \ref{L:step0}, $\theta^{(k)}$ has the role of $\theta$ (provided that $I(\theta^{(k)})>0$)
and $\theta^{(k+1)}$ has the role of $\widetilde \theta$.

With each recursion
$$
\ell(\theta^{(k)} \mid \mathbf x)<\ell(\theta^{(k+1)} \mid \mathbf x)
$$
and either $S(\theta^{(k+1)})$ has
fewer members than $S(\theta^{(k)})$, or
$$
S(\theta^{(k+1)})=S(\theta^{(k)}) \quad \text{and} \quad I(\theta^{(k+1)})<I(\theta^{(k)})
$$
Because $S(\theta^{(0)})$ is finite, 
we can reduce the
size of $S(\theta^{(k)})$ at most finitely many times.  Let $K$ be the last time
that $S(\theta^{(K)})$ is smaller than $S(\theta^{(K-1)})$.  Then,
after $M=I(\theta^{(K)})$ more steps of the recursion, $I(\theta^{(K+M)})=0$.
\end{proof}

\begin{theorem}  \label{T:exist} There is a non-decreasing response estimate $\theta \in \Theta^m$
such that $I(\theta)=0$, $\ell(\theta \mid \mathbf x)>0$  and, for all non-decreasing $\lambda \in \Theta^m$,
$$
\ell(\lambda \mid \mathbf x) \le \ell(\theta \mid \mathbf x)
$$
\end{theorem}
\begin{proof}  By Lemma \ref{L:step1}, for any non-decreasing $\lambda \in \Theta^m$, there is some non-decreasing $\theta \in \Theta^m$ with $I(\theta)=0$ such that
$$
\ell(\lambda \mid \mathbf x) \le \ell(\theta \mid \mathbf x)
$$

Note that, when $I(\theta)=0$, the partition $S(\theta)$ induced by $\theta$ determines
$\theta$:
$$
i \in A \in S(\theta) \quad \rightarrow \quad \theta_i = \mu(A)
$$

Because $[1 \hdots m]$ is finite, there are finitely many distinct partitions of it.
With the condition $I(\theta)=0$ imposed, that gives at most finitely many members of $\Theta^m$
that can play the role of $\tilde \theta$ in Lemma \ref{L:step1}.  Among these finitely many possibilities for $\widetilde \theta$ in Lemma \ref{L:step1}, choose some $\widehat \theta$ with $\ell(\widehat \theta \mid \mathbf x)$ being largest.

In particular, consider the $\lambda \in \Theta^m$ that is constant with components equal to  $\overline x=\mu[1 \hdots m]$.  Because all $x_{i,j}$ are observable, by Lemma \ref{L:begin}
we have 
$$\ell(\lambda \mid \mathbf x)=\ell(\overline x \mid \mathbf x)>0$$  It follows that
$$
\ell(\widehat \theta \mid \mathbf x) \ge \ell(\lambda \mid \mathbf x)>0
$$
\end{proof}

\section{Algorithm and Uniqueness}

\begin{lemma} \label{L:splitting} Suppose that all $x_{i,j} \in \mathcal D$.  Let $\theta$ satisfy the conclusions of Theorem \ref{T:exist}.  Let $\tau_1 < \tau_2 < \hdots <\tau_s$ be
a complete listing of the distinct components of $\theta$.  For $r \in [1 \hdots s]$, set
$$
A_r=\setof{i \in [1 \hdots m]}{\theta_i=\tau_r}
$$
Then
\begin{enumerate}
\item
There are integers $a_r \le b_r$ in $[1 \hdots m]$ such that $A_r=[a_r \hdots b_r]$.
\item
For all $t \in[(a_r+1)\hdots b_r]$,
$$
\mu[t \hdots b_r] \le \mu(A_r)
$$
\end{enumerate}
\end{lemma}
\begin{proof}  Because $\theta$ is non-decreasing, Item (i) is immediate.

We will prove Item(ii) by contradiction.  Suppose there is $r$ and $t \in[(a_r+1) \hdots b_r]$ such that
$$
\mu[t \hdots b_r]>\mu(A_r)
$$
Because $I(\theta)=0$, we have $\mu(A_r)=\tau_r$.

For some real $\eta$ to be chosen soon, we define $\widetilde \theta$
as follows: for $i \in [1 \hdots m]$ let
$$
\widetilde{\theta_i}=\begin{cases}
\theta_i &\text{if $i<t$ or $i>b_r$} \\
\eta & \text{if $i \in [t \hdots b_r]$}
\end{cases}
$$

If $r=s$, let $\eta=\mu[t \hdots b_r]$.  If $r<s$, select
$$\eta \in (\tau_r,\min\{\mu[t \hdots b_r], \tau_{r+1}\})
$$
Please note that this selection makes $\widetilde \theta$ be
non-decreasing.  Also, $\eta \in \Theta$:
\begin{enumerate}
\item
Because every $x_{i,j}$ is an observable real number, $\eta=\mu[t \hdots b_r]$ puts $\eta \in \Theta$ by Lemma 1.  
\item
Because $\Theta$ is a real interval and both
$\tau_r$ and $\tau_{r+1}$ are in $\Theta$, $\eta \in (\tau_r,\tau_{r+1})$ puts $\eta \in \Theta$.
\end{enumerate}
Therefore $\widetilde \theta$ is 
a non-decreasing member of $\Theta^m$.

Let
$$
C=\left(\prod_{i=1}^{t-1} \prod_{j=1}^{n_i} f(x_{i,j} \mid \theta_i) \right) \cdot
\left(\prod_{i=b_r+1}^m \prod_{j=1}^{n_i} f(x_{i,j} \mid \theta_i)\right)
$$
(with empty products set equal to $1$).   Since $C \ge 0$ and a factor of
$\ell(\theta \mid \mathbf x)>0$, we have $C>0$.

By Assumption 3, with  $\tau_r < \eta \le \mu[t \hdots b_r]$,
$$
\begin{aligned}
\prod_{i=t}^{b_r} \prod_{j=1}^{n_i} f(x_{i,j} \mid \eta)
&> \prod_{i=t}^{b_r} \prod_{j=1}^{n_i} f(x_{i,j} \mid \tau_r)\\
&=\prod_{i=t}^{b_r} \prod_{j=1}^{n_i} f(x_{i,j} \mid\theta_i)
\end{aligned}
$$
Therefore,
$$
\ell(\widetilde \theta \mid \mathbf x)=
C \cdot \prod_{i=t}^{b_r} \prod_{j=1}^{n_i} f(x_{i,j} \mid \eta)
> C \cdot \prod_{i=t}^{b_r} \prod_{j=1}^{n_i} f(x_{i,j} \mid \theta_i)
=\ell(\theta \mid \mathbf x)
$$
This contradicts the likelihood maximizing property of $\theta$ in Theorem \ref{T:exist}.
\end{proof}

\begin{theorem}\label{T:algo}Suppose that all $x_{i,j} \in \mathcal D$.  Let $\theta$ satisfy the conclusions of Theorem \ref{T:exist}.  Let $\tau_1 < \tau_2 < \hdots <\tau_s$ be
a complete listing of the distinct components of $\theta$.  For $r \in [1 \hdots s]$, set
$$
A_r=\setof{i \in [1 \hdots m]}{\theta_i=\tau_r}
$$
Then
\begin{enumerate}
\item
There are integers $a_r \le b_r$ in $[1 \hdots m]$ such that $A_r=[a_r \hdots b_r]$.
\item
Let 
$$
\kappa_r=\min \setof{\mu[a_r \hdots k]}{k \in [a_r \hdots m]}
$$
and set
$$
t_r=\max\setof{k \in [a_r \hdots m]}{\mu[a_r \hdots k]=\kappa_r}
$$
Then $t_r=b_r$.
\end{enumerate}
\end{theorem}

\begin{proof}  Item (i) is the same as Item(i) in Lemma \ref{L:splitting}, and is repeated here to establish the notation.

Because $I(\theta)=0$, we have $\tau_r=\mu(A_r)$ for all $r$.

Suppose first there is some $r$ such that $t_r=b_q$ for some $q>r$.
Then $\mu[a_r \hdots t_r]$ is a convex combination with positive coefficients of $\tau_h$ for $h \in [r \hdots q]$.  For $h>r$, we have $\tau_h>\tau_r$.  By the definition of
$\kappa_r$, we have
$$
\tau_r=\mu(A_r)=\mu[a_r \hdots b_r] \ge \kappa_r
$$
Therefore, $\kappa_r>\kappa_r$.  This contradiction proves that $t_r \ne b_q$ for
all $q \in [r+1 \hdots s]$.

Next suppose that $t_r \in [a_q \hdots (b_q-1)]$ for some $q \in [r \hdots s]$.  By the definition of $\kappa_r$ and of $t_r$, 
\begin{itemize}
\item
$\mu[a_r \hdots b_q]>\kappa_r$
\item
When $q>r$, 
$$
h \in [r \hdots (q-1)] \quad \Rightarrow \quad \mu[a_r \hdots b_h] \ge \kappa_r
$$
In particular, $\mu[a_r \hdots b_{q-1}] \ge \kappa_r$.
\end{itemize}

We now show that $\mu[a_q \hdots t_r] \le \kappa_r$.  If $q=r$ this is immediate
from the definitions of $\kappa_r$ and $t_r$.
Suppose that $q>r$.  There are positive integers $e$ and $f$ such that
$$
\mu[a_r \hdots t_r]=\dfrac{e \mu[a_r \hdots b_{q-1}] +f \mu[a_q \hdots t_r]}{e+f}
$$
Therefore
$$
\begin{aligned}
\mu[a_q \hdots t_r]&=\dfrac{(e+f) \mu[a_r \hdots t_r]-e \mu[a_r \hdots b_{q-1}]}{f} \\
&=\dfrac{(e+f) \kappa_r-e \mu[a_r \hdots b_{q-1}]}{f}\\
&=\kappa_r +(e/f) \cdot \left\{\kappa_r-\mu[a_r \hdots b_{q-1}] \right\} \\
&\le \kappa_r
\end{aligned}
$$
because $\mu[a_r \hdots b_{q-1}] \ge \kappa_r$.

Next we show that $\mu[(t_r+1) \hdots b_q]>\kappa_r$.  There are positive
integers $e$ and $f$ such that 
$$
\mu[a_r \hdots b_q]=\dfrac{e \mu[a_r \hdots t_r] +f  \mu[(t_r+1) \hdots b_q]}{e+f}
$$
Consequently
$$
\begin{aligned}
\mu[(t_r+1)\hdots b_q]&=\dfrac{(e+f)\mu[a_r \hdots b_q]-e \mu[a_r \hdots t_r]}{f} \\
&=\dfrac{(e+f)\mu[a_r \hdots b_q]-e \kappa_r}{f}\\
&=\mu[a_r \hdots b_q]+(e/f) \cdot\left \{\mu[a_r \hdots b_q]-\kappa_r \right\}\\
&>\mu[a_r \hdots b_q] > \kappa_r
\end{aligned}
$$
because $\mu[a_r \hdots b_q]>\kappa_r$.

Third, we argue that $\mu[(t_r+1) \hdots b_q]>\tau_q$.  We've shown already
that
$$
\mu[a_q \hdots t_r] \le \kappa_r < \mu[(t_r+1) \hdots b_q]
$$
There are positive integers $e$ and $f$
such that
$$
\tau_q=\mu[a_q \hdots b_q]=\dfrac{e \mu[a_q \hdots t_r]+f\mu[(t_r+1) \hdots b_q]}{e+f}
$$
It follows that
$$
\tau_q < \dfrac{e \mu[(t_r+1) \hdots b_q] +f \mu[(t_r+1) \hdots b_q]}{e+f} =
\mu[(t_r+1) \hdots b_q]
$$
However, this contradicts Lemma \ref{L:splitting} (since $\theta$ satisfies the conclusions of Theorem \ref{T:exist}).  

The only possibility left for $t_r$ is to be equal to $b_r$ as desired.
\end{proof}

\begin{corollary}  Suppose that all $x_{i,j} \in \mathcal D$.  There is a unique $\theta \in \Theta^m$ that is non-decreasing and, for all non-decreasing $\lambda \in \Theta^m$,
\begin{equation}
\ell(\lambda \mid \mathbf x) \le \ell(\theta \mid \mathbf x)
\label{E:eqn}
\end{equation}
where $\ell$ is the compound likelihood function.
\end{corollary} 
\begin{proof}  By Theorem \ref{T:exist}, there is some $\theta$ such that satisfies the conclusions of that theorem.   In particular, $\theta$ satisfies Equation \ref{E:eqn}
for all non-decreasing $\lambda \in \Theta^m$.

Suppose that $\widehat \theta$ also satisfies Equation \ref{E:eqn} for all non-decreasing
$\lambda \in \Theta^m$.  By Lemma \ref{L:step1}, if $I(\widehat \theta)>0$,
there would be some non-decreasing $\widetilde \theta \in \Theta^m$ such that
$$
\ell(\widetilde \theta \mid \mathbf x)> \ell(\widehat\theta \mid \mathbf x)
$$
So we must have $I(\widehat \theta)=0$.

 Because each $x_{i,j}$ is observable, their arithmetic mean $\mu[1 \hdots m]$ is in $\Theta$.  Let $\lambda \in \Theta^m$ which has the constant component $\mu([1 \hdots m])$.  By Assumption 3,
$$
\ell(\lambda \mid \mathbf x)>0
$$
Therefore, $\ell(\widehat \theta \mid \mathbf x)>0$.

So $\widehat \theta$ satisfies the conclusions of Theorem \ref{T:exist}.

Note that Theorem \ref{T:algo} specifies $S(\theta)=S(\widehat \theta)$ uniquely.
First, the theorem determines $A_1=[1 \hdots b_1]$ as $b_1$ must equal $t_1$.
Once $b_r$ is determined, if $b_r<m$ the theorem then determines $b_{r+1}=t_{r+1}$.  Of course $a_{r+1}=b_r+1$, and thus $A_{r+1}$ is specified.

However, since $I(\theta)=I(\widehat \theta)=0$, we have for $i \in A_r$
$$
\theta_i =\widehat{\theta_i} = \mu(A_r)
$$
Thus $\theta=\widehat{\theta}$.
\end{proof}

\section{Examples of the Formal Context}

Throughout this section, let $T$ be a positive integer, $\mathbf y:[1 \hdots T] \rightarrow \mathcal D$, and we think of each $y_i$ as an observation of a random variable $Y_i$, with
$\{Y_i\}_{i=1}^T$ independent with PDFs from the given family.

\begin{example}  Let $\Theta=[0,1]$ and $\mathcal F$ be the family of Bernoulli random variables (here including two constant random variables).
Given $\theta \in \Theta$, the PDF with parameter $\theta$ is defined as follows:
$$
f(x\mid\theta) =\begin{cases}
\theta &\text{if $x=1$} \\
1-\theta & \text{if $x=0$} \\
0 &\text{if $x \in \mathbb R\backslash \{0,1\}$}
\end{cases}
$$
\end{example}

Here $\mathcal D=\{0,1\}$.  We show that Assumption 3 holds.

Let $\mathbf y$ have $r$ ones and $s$ zeros.  Thus $r+s=T$ and $\overline y=r/(r+s)$.
For any $\lambda \in [0,1]$,
$$
h(\lambda):=L(\mathbf \lambda \mid\mathbf y)=\lambda^r (1-\lambda)^s
$$

{\bf Case 1:}  Suppose $r=0$.  So $s=T>0$.  Then
$$
h'(\lambda) = s (1-\lambda)^{s-1} (-1)
$$
For $\lambda \in [0,1)$, this derivative is negative.  Because $h(\lambda)$ is continuous in
$\lambda$, the function $h$ is strictly decreasing on $[0,1]$.
Since $\mathbf y$ has all zeros, $\overline{y}=0$.  So Assumption 3 holds:  $h$ is strictly decreasing on $[0,1]=[\overline y, 1]$
and strictly increasing (trivially) on $[0,0]$.

{\bf Case 2:}  Suppose $s=0$.  Then $r=T>0$ and
$$
h'(\lambda)=r \lambda^{r-1}
$$
For $\lambda \in (0,1]$, this derivative is positive. Because $h(\lambda)$ is continuous in
$\lambda$, the function $h$ is strictly increasing on $[0,1]$.
Since $\mathbf y$ has all ones, $\overline{y}=1$.  Consequently, Assumption 3 holds:  $h$ is strictly increasing
on $[0,\overline y]=[0,1]$ and $h$ is strictly decreasing (trivially) on $[1,1]$.

{\bf Case 3:}  Suppose $r>0$ and $s>0$.  Then
$$
\begin{aligned}
h'(\lambda) & = r \lambda^{r-1} (1-\lambda)^s + \lambda^r \cdot s (1-\lambda)^{s-1} (-1) \\
&=\lambda^{r-1} (1-\lambda)^{s-1} \left[r (1-\lambda)-s \lambda \right] \\
&=\lambda^{r-1} (1-\lambda)^{s-1} (r+s)\left[\dfrac{r}{r+s}-\lambda\right] \\
\end{aligned}
$$
For $\lambda \in (0,1)$, this derivative is positive for $\lambda<r/(r+s)$ and negative for $\lambda>r/(r+s)$.
Because $h(\lambda)$ is continuous in $\lambda$, it follows that $h$ is strictly increasing on $[0,r/(r+s)]$
and strictly decreasing on $[r/(r+s),1]$.  Since $\overline y=r/(r+s)$, Assumption 3 holds.

\begin{example}  Let $\Theta=[0,\infty)$ and $\mathcal F$ be the family of Poisson random variables (here including one constant random variable).
Given $\theta \in \Theta$, the PDF with parameter $\theta$ is defined as follows:
$$
f(x\mid\theta) =\begin{cases}
\dfrac{e^{-\theta}\theta^x}{x!}, & \text{if $x \ge 0$ and an integer} \\
\\
0 &\text{otherwise}
\end{cases}
$$
\end{example}

Here $\mathcal D$ is the set of non-negative integers.  We show that Assumption 3 holds.

For any $\lambda \in [0,\infty)$, let
$$
h(\lambda):=L(\lambda\mid\mathbf y)=\prod_{i=1}^T \left(e^{-\lambda} \dfrac {\lambda^{y_i}}{y_i !} \right)=K e^{-T\lambda} \lambda^{T \overline{y}}
$$
where $K>0$ is a factor that does not depend on $\lambda$.

{\bf Case 1.}  Suppose that $\overline{y}=0$.  Then $h(\lambda)=Ke^{-T \lambda}$ and hence
$$
h'(\lambda)=-KT e^{-T\lambda}
$$
Note that $h'$ is negative for all $\lambda \in \Theta$.  Thus $h$ is strictly decreasing on $[0, \infty)=[\overline y,\infty)$
and strictly increasing (trivially) on $[0,0]$.

{\bf Case 2.}  Suppose that $\overline{y}>0$.
Then
$$
\begin{aligned}
h'(\lambda)&=(-KT) e^{-T\lambda} \lambda^{T \overline{y}} + Ke^{-T \lambda} \cdot T \overline{y} \cdot \lambda^{T \overline{y}-1}\\
&=KTe^{-T \lambda} \lambda^{T \overline{y}-1} \left[-\lambda + \overline{y} \right]
\end{aligned}
$$
Note that $h'$ is positive on $(0,\overline{y})$ and $h'$ is negative on $(\overline{y}, \infty)$.  Since  $\overline{y} >0$, we have
$h$ continuous on $[0,\infty)=\Theta$.  Therefore, $h$ is strictly increasing
on $[0,\overline{y}]$ and $h$ is strictly decreasing on $[\overline{y},\infty)$.

\begin{example}
Let $\Theta=[0, \infty)$ and $\mathcal F$ be the family of geometric random variables (here including one constant random variable), parameterized by
their means.
Given $\theta \in \Theta$, set $p=1/(1+\theta)$.  For $\theta >0$, the PDF with parameter $\theta$ is defined as follows:
$$
f(x\mid\theta) =\begin{cases}
(1-p)^x p, & \text{if $x \ge 0$ and an integer} \\
\\
0 &\text{otherwise}
\end{cases}
$$
For $\theta=0$, and thus $p=1$, let $f(x\mid 0)=0$ for all $x$ except that $f(0\mid 0)=1$.
\end{example}

Here $\mathcal D$ is the set of non-negative integers.

For any $\lambda \in [0,\infty)$, set $p=1/(1+\lambda)$.
Then for $\lambda>0$
$$
h(\lambda):=L(\lambda \mid\mathbf y)=\prod_{i=1}^T \left((1-p)^{y_i} p\right)
=p^T (1-p)^{T \overline{y}}
$$
For $\lambda=0$, $h(0)=1$ if $\mathbf y$ is a vector of zeros and $0$ otherwise.

Please note that $h$ is continuous on $[0, \infty)$:
\begin{itemize}
\item
Suppose $\overline{y}=0$.  Then $\mathbf y$ is a vector of zeros and $h(0)=1$.  For $\lambda>0$, we have
$$
h(\lambda)=p^T(1-p)^{T \cdot 0}=p^T=(1+\lambda)^{-T}
$$
Clearly $h$ is continuous on $(0,\infty)$; it is also continuous at $0$ because $\lim_{\lambda \downarrow 0} h(\lambda)=1=h(0)$.
\item
Suppose $\overline{y}>0$.  Then $\mathbf y$ has at least one non-zero component and thus $h(0)=0$.  For $\lambda >0$,
$$
h(\lambda)=\lambda^{T \overline{y}} (1+\lambda)^{-T \overline{y}-T}
$$
Clearly $h$ is continuous on $(0,\infty)$.  Since $T \overline{y}>0$, we have $\lim_{\lambda \downarrow 0} h(\lambda)=0=h(0)$.
Thus $h$ is continuous at $0$ as well.
\end{itemize}

{\bf Case 1.}  Suppose that $\overline{y}=0$.  Then $\mathbf y$ is a vector of zeros and $h(0)=1$.  For $\lambda>0$,
$$
h(\lambda)=p^T=(1+\lambda)^{-T} \quad \text{and thus}\quad h'(\lambda)=(-T)(1+\lambda)^{-T-1}<0
$$
Since $h$ is continuous on $[0,\infty)$ and has a negative derivative on $(0,\infty)$, we know that $h$ is strictly decreasing on $[0,\infty)$.
Thus $h$ is strictly decreasing on $[0,\infty)=[\overline y,\infty)$ and strictly increasing
(trivially) on $[0,0]$.

{\bf Case 2.}  Suppose that $\overline{y}>0$.  For $\lambda>0$,
$$
\begin{aligned}
h'(\lambda)
&= \dfrac{dh}{dp} \cdot \dfrac{dp}{d\lambda} \\
&=\left\{T p^{T-1} (1-p)^{T \overline{y}} + p^T \cdot T \overline{y} \cdot (1-p)^{T \overline{y}-1}(-1)\right \}\\
&\quad \quad  \cdot (-1)(1+\lambda)^{-2}\\
&=-(1+\lambda)^{-2} T p^{T-1} (1-p)^{T \overline{y}-1} \left[ (1-p)-\overline{y}p \right] \\
&=(1+\lambda)^{-2}T p^{T-1} (1-p)^{T \overline{y}-1} \left[\overline{y}p-(1-p) \right]\\
&=(1+\lambda)^{-2}T p^{T-1} (1-p)^{T \overline{y}-1} \left[ \dfrac{1+\overline{y}}{1+\lambda}-1 \right] \\
&=(1+\lambda)^{-2}T p^{T-1} (1-p)^{T \overline{y}-1} \left[ \dfrac{\overline{y}-\lambda}{(1+\lambda)(1+\overline{y})} \right]
\end{aligned}
$$
For $\lambda >0$, all the factors in the previous line are positive except for the factor $\overline{y}-\lambda$.
Thus for $\lambda \in (0,\overline{y})$, we have $h'(\lambda)>0$ and for $\lambda \in (\overline{y},\infty)$ we have $h'(\lambda)<0$.
Since $h$ is continuous on $[0, \infty)$, we have $h$ strictly increasing on $[0, \overline{y}]$ and strictly decreasing on $[\overline{y},\infty)$.

\begin{example}
Let $\Theta=\mathbb R$ and, for a fixed $\sigma>0$, let $\mathcal F$ be the family of normal random variables with standard deviation $\sigma$.
Given $\theta \in \Theta$, the PDF with parameter $\theta$ is defined as follows: for all real $x$
$$
f(x\mid\theta) = \dfrac{1}{\sigma \sqrt{2 \pi}} \exp\left\{\dfrac{-(x-\theta)^2}{2 \sigma^2}\right\}
$$
\end{example}

Here $\mathcal D =\mathbb R$.  We'll argue that Assumption 3 holds.

For all real $\lambda$, let
$$
\begin{aligned}
h(\lambda)&:=L(\lambda \mid\mathbf y)=\prod_{i=1}^T \left(\dfrac{1}{\sigma \sqrt{2 \pi}} \exp\left\{\dfrac{-(y_i-\theta)^2}{2 \sigma^2}\right\}\right)\\
&=K \exp \left\{\dfrac{-\sum_{i=1}^T (y_i-\lambda)^2}{2 \sigma^2}\right\}
\end{aligned}
$$
where $K$ is a positive factor that does not depend on $\lambda$.  Then
$$
\begin{aligned}
h'(\lambda)&= K \exp \left\{\dfrac{-\sum_{i=1}^T (y_i-\lambda)^2}{2 \sigma^2}\right\} \cdot \left(\dfrac{\sum_{i=1}^T 2(y_i-\lambda)}{2 \sigma^2} \right)
\\
&=\dfrac{K}{\sigma^2} \cdot \exp \left\{\dfrac{-\sum_{i=1}^T (y_i-\lambda)^2}{2 \sigma^2}\right\} \cdot \left(-T \lambda+\sum_{i=1}^T y_i \right) \\
&=\dfrac{TK}{\sigma^2} \cdot \exp \left\{\dfrac{-\sum_{i=1}^T (y_i-\lambda)^2}{2 \sigma^2}\right\} \cdot \left(-\lambda + \overline{y} \right)
\end{aligned}
$$
All factors immediately above for $h'(\lambda)$, except the last one,  are positive for all $\lambda$.
So $h'$ is positive for $\lambda$ in $(-\infty, \overline y)$ and
$h'$ is negative for $\lambda$ in $(\overline y,\infty)$.  Since $h$ is continuous on $\mathbf R$, we have $h$ is strictly increasing on $(-\infty, \overline{y}]$ and strictly decreasing
on $[\overline{y},\infty)$.

\begin{example}
Let $\Theta=(0, \infty)$ and $\mathcal F$ be the family of exponential random variables parameterized by their means.
Given $\theta \in \Theta$, set $\tau=1/\theta$.  In terms of $\tau$, the PDF with parameter $\theta$ is defined as follows:
$$
f(x\mid\theta) =\begin{cases}
\tau e^{-\tau x}, & \text{for  $x > 0$} \\
\\
0 &\text{for $x \le 0$}
\end{cases}
$$
\end{example}
Here $\mathcal D=(0,\infty)$.  We'll argue that Assumption 3 holds.

Note that, since $\mathcal D$ has only positive numbers, $\overline y>0$.

For any real $\theta>0$ and with $\tau=1/\theta$, let
$$
h(\theta):=L(\theta \mid\mathbf y)=\prod_{i=1}^T \left(\tau e^{-\tau y_i}\right)
=\tau^T e^{-\tau T \overline y}
$$
Then
$$
\begin{aligned}
h'(\theta)&=\dfrac{dh}{d\tau} \cdot \dfrac{d \tau}{d\theta} \\
&= \left\{ T \tau^{T-1} e^{-\tau T \overline y} + \tau^T \cdot (-T \overline y) \cdot e^{-\tau T \overline y} \right\}
\cdot \left(-\theta^{-2} \right) \\
&=\dfrac{-T \tau^{T-1} e^{-\tau T \overline y} }{\theta^2} \left\{1-\tau \overline y \right\} \\
&=\dfrac{T \tau^{T-1} e^{-\tau T \overline y}}{\theta^2} \left\{ \dfrac{\overline y - \theta}{\theta} \right\}
\end{aligned}
$$
It follows that $h'(\theta)>0$ if if $\theta<\overline y$ and $h'(\theta)<0$ if $\theta>\overline y$.  Because $h(\theta)$ is continuous on $(0,\infty)$, we have 
$h$ strictly increasing on $(0,\overline y]$ and strictly decreasing on $[\overline y,\infty)$. 

\newpage

\section{Data Appendix}

Table \ref{tab:data} has the no-show counts, tabulated by SAT-R scores.

\begin{table}
\centering
\begin{tabular}{cccc}
\text{SAT-R Scores} &\text{Total Count} &\text{No-Show Count} \\
330 &1 &0&\\
390 & 2&0&\\
400 & 1 &0&\\
410& 2&0&\\
420& 5&0&\\
430&4&0&\\
440&4&1&\\
450&3&2&\\
460&2&0&\\
470&8&1&\\
480&11&3&\\
490&9&0&\\
500&4&1 &\\
510&11&0&\\
520&9&0&\\
530&8&1&\\
540&11&4&\\
550&6&1&\\
560&5&0&\\
570&6&0&\\
580&7&0&\\
590&5&1&\\
600&3&1&\\
610&5&3&\\
620&4&2&\\
630&1&0&\\
640&7&2&\\
650&1&0&\\
660&1&0&\\
680&1&1&\\
690&1&1&\\
700&1&0&\\
710&1&0&\\
750&1&0&\\
800&1&1&\\
\end{tabular}
\caption{No-Show Counts for SAT-R Levels}
\label{tab:data}
\end{table}

\end{document}